%% file: agt-1-16.tex
\documentclass{gtart}
\input agtout

\lognumber{16}
\volumenumber{1}
\volumeyear{2001}
\papernumber{16}
\pagenumbers{321}{347}
\received{14 December 2000}
\revised{21 May 2001}
\accepted{25 May 2001}
\published{31 May 2001}

\usepackage{amsmath,amssymb} 
\usepackage{graphics,labelfig}

\newtheorem{theorem}{Theorem}[subsection]

\newtheorem{proposition}[theorem]{Proposition}
\newtheorem{lemma}[theorem]{Lemma}
\newtheorem{corollary}[theorem]{Corollary}
\newtheorem{scholium}[theorem]{Scholium}
\newtheorem{remark}[theorem]{Remark} 
\newtheorem {conjecture}[theorem]{Conjecture}
\newtheorem {example}[theorem]{Example}

\newenvironment {acknowledgement}{\noindent {\bf
    Acknowledgment}}{\medskip}

\newcommand{\Z} {\mathbb{Z}}
\newcommand{\R} {\mathbb{R}} 
\newcommand{\s}{{\sigma}} 
\newcommand{\tn} [1] {\| #1 \|_T} 
\newcommand{\an} [1] {\| #1 \|_A}

\newcommand{\BeT} {\beta_t}   
\newcommand{\br} {\mathfrak{b}}     
\newcommand{\id} {\mbox{id}} 
\newcommand{\HOMFLY} {\mbox{HOMFLY}}
\newcommand{\pos} {\mbox{pos}}
\newcommand{\negat} {\mbox{neg}}

\begin{document} 

\title[On McMullen's and other inequalities]{On McMullen's and other inequalities\\for the Thurston norm of 
link complements}

\author {Oliver T. Dasbach\\Brian S. Mangum}
\address {University of California, Riverside, Department of
  Mathematics\\Riverside, CA 92521 - 0135, USA}
\email{kasten@math.ucr.edu}
\url{http://www.math.ucr.edu/\char'176kasten}

\secondaddress {Barnard College/Columbia University,
Department of Mathematics\\New York, NY 10027, USA}
\secondemail{mangum@math.columbia.edu}
\asciiemail{kasten@math.ucr.edu, mangum@math.columbia.edu}
\asciiaddress{University of California, Riverside, Department of
  Mathematics\\Riverside, CA 92521 - 0135, USA\\
Barnard College/Columbia University,
Department of Mathematics\\New York, NY 10027, USA}


\begin{abstract} 
  In a recent paper, McMullen showed an inequality between the
  Thurston norm and the Alexander norm of a $3$-manifold.
  This generalizes the well-known fact that twice the genus
  of a knot is bounded from below by the degree of the Alexander
  polynomial.

  We extend the Bennequin inequality for links to an
  inequality for all points of the Thurston norm, if the manifold
  is a link complement. We compare these two inequalities on
  two classes of closed braids.

  In an additional section we discuss a conjectured inequality due
  to Morton for certain points of the Thurston norm.
  We prove Morton's conjecture for closed $3$-braids.
\end {abstract}


\asciikeywords{Thurston norm, Alexander norm, multivariable Alexander 
polynomial, fibred links, positive braids, Bennequin's inequality, 
Bennequin surface, Morton's conjecture}

\maketitle

Recently, the multivariable Alexander polynomial experienced a
fulminant\break comeback on the mathematical stage. Connections to
Donaldson-Seiberg-Wit\-ten theory and also its meaning in terms of
Vassiliev knot invariants became known.  In this paper we focus our
attention on the recent inequality of McMullen \cite
{McMullen:ThurstonAlexanderNorm} that compares the Alexander norm of a
$3$-manifold $M$ with its Thurs\-ton norm. Both norms are norms on
$H^1(M;\R)$. If the manifold is the complement  of an oriented link
$L$ in $S^3$ or, more generally, in an integral homology sphere 
then the meridians of $L$ give a natural base for $H^1(M,\R)$. 
The Thurston norm at the point $(1,\dots,1)$ is determined by a minimal
spanning surface of the link.
McMullen shows that if the first Betti number of the
manifold is at least $2$ then the Thurston norm is greater than or equal
to
the Alexander norm on all points. Moreover, he compares the value on 
all complements of links with less than $9$ crossings.  
McMullen's results extend the well-known fact that twice the genus of
a knot is greater than or equal to the span of the Alexander
polynomial.

On the other hand, another inequality for the Thurston norm of link
complements, evaluated at certain points, is sometimes quite useful:
The Bennequin inequality.
We will extend this inequality to an inequality for all points of
the Thurston norm. 

To show the differences, we will compare these two inequalities on two 
classes of links.

\begin{enumerate}
\item The first class is the class of closed homogeneous braids,
i.e.\ braids in which the standard generator $\s_i$ appears
at least once for every $i$, and all exponents of $\s_i$ have the same
sign.
By a result of Stallings \cite{Stallings:fibred_links} 
these links are known to be fibred and the fibre surface is the one 
that one gets by applying Seifert's algorithm.

This implies that the Thurston and the Alexander norm of the complements
of these links agree at the points given by the fibration. 
With some additional arguments we will give the equality for all points
$C=(C_1, \dots, C_r) \in H^1(M,\R)$ with $C_i \geq 0$.
In contrast to this result, Dunfield \cite{Dunfield:McMullen_estimate}
shows that the equality does not necessarily hold at all points for
three-manifolds that fibre over the circle.

The Bennequin inequality is in general not helpful for this class of links.
We will give, however, a simple proof of a recent result \cite{Kanda}
saying that the difference between the Bennequin number and the
Thurston norm can be arbitrarily large.

\item The second class is the class of closed band
positive braids. The band generators, a new set of generators of the 
braid groups,
have recently appeared quite useful \cite{BKL}. 
One defines $a_{i,j}$ to be the
braid 
$$\s_i^{-1} \cdots \s_{j-2}^{-1} \s_{j-1} \s_{j-2} \cdots \s_i.$$
With these band generators, band positive braids are braids where
all $a_{i,j}$ appear only with positive exponents.

Ko and Lee \cite{KL:genera} showed that a band Seifert algorithm 
for closed band positive braids on four strands yields a minimal spanning
surface. By computing the Alexander polynomial of 
an example they point out that furthermore these band positive braids 
are not necessary fibred. From this example and others, it 
follows that McMullen's inequality is not sharp for this
class of links.

The Bennequin inequality gives the 
Thurston norm at the point $C=(1,\dots,1)$
for closed band positive braids of arbitrary braid index.
In particular the result of Ko and Lee holds for an arbitrary number of 
strands. It follows that there is 
a Bennequin surface, in the sense of \cite{BirmanMenasco:Bennequin_surface},
 for every closed band positive braid.  For special classes of band positive
braids, the generalization of Bennequin's inequality will give us the
Thurston norm at all points $C=(C_1,\dots,C_r), C_i \geq 0$.
\end{enumerate}

In addition to these two inequalities there
is a conjectured inequality, due to Morton, for points
$C=(C_1,\dots,C_r), \vert C_i \vert=1,$ of the 
Thuston norm that comes from the HOMFLY polynomial. 
This inequality would sharpen the Bennequin inequality.

To prove it, however, seems to be a hard task. 
It is known that this conjecture is true for all closed homogeneous braids. 
We will give a new and different  proof of this result, and
we will prove that it is also true for the class of closed $3$-braids.

Our proof of this latter result uses powerful tools, as a beautiful theorem
of Scharlemann and Thompson, based on Gabai's work.

Added in proof: It was well-known, that the span of the
Alexander polynomial does not give a lower bound for the slice-genus
of a knot. The slice-genus is the minimal genus of an oriented smooth embedded
surface $F$ in $D^4$ with boundary the knot. 
Recently, Ferrand \cite{Ferrand:Legendrian} gave
an example of a knot where a ``slice Morton conjecture''
does not hold. The slice Bennequin inequality, however, was proven in
\cite{Rudolph:SliceBennequin}.
 
\bigskip
\begin{acknowledgement}
  The work to this paper was mainly done while the first author was
  visiting Columbia University in Summer 1999.
  The authors would like to express their deep gratitude to Joan Birman 
  for many helpful discussions and suggestions that influenced this
  paper. We also appreciate very much the efforts of an anonymous 
  referee to make the paper more readable. 

 Furthermore, the first author thanks Hugh Morton, Guennadi Noskov 
 and Arkady Vaintrob for their helpful remarks.  

 The work of the first author was partially supported
 by the Deutsche Forsch\-ungsgemeinschaft (DFG).
\end{acknowledgement}

\section{Setting the scene}
\label{setting_the_scene}

\subsection{A generalization of Seifert's algorithm}
\label{generalized Seifert algorithm}

\begin{figure}[ht]
\begin{center} 
\parbox{2 cm}{\scalebox{0.37}{\includegraphics{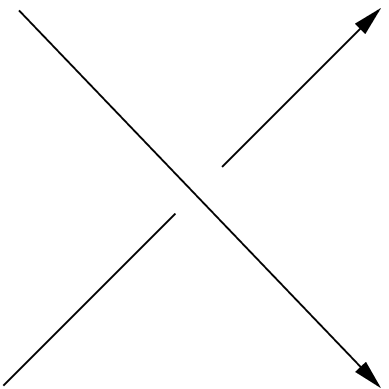}}}
or
\parbox{2cm} {\scalebox{0.37}{\includegraphics{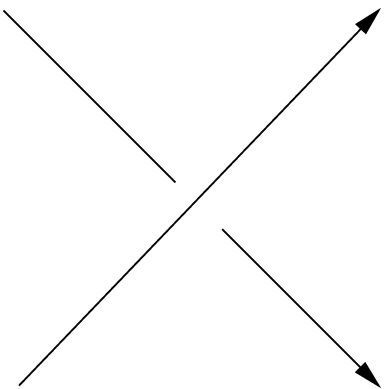}}}
$\longrightarrow$
\parbox{2cm}{\scalebox{0.5}{\includegraphics{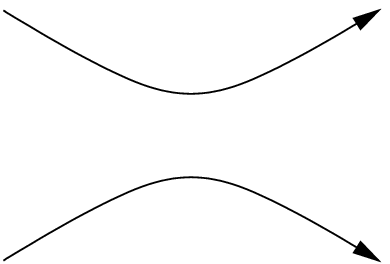}}}
\end{center}
\caption{\label{smoothing} Smoothing a crossing}
\end{figure}

By a classical result of Seifert, every oriented link is the boundary of
an oriented surface in $S^3$. His algorithm to construct such a surface is
by smoothing every crossing in a diagram of the link (see Figure
\ref{smoothing}). This gives a
collection of circles, called Seifert circles. Now bands are added
to join the Seifert circles meeting at a crossing. The surface so 
constructed has Euler characteristic $(s - c)$ where $c$ is the number of
crossings and $s$ is the number of Seifert circles.

One can generalize Seifert's algorithm to find a punctured oriented surface 
spanning a sublink of a link in the complement of the other components:
For the sublink construct a surface as described above. Whenever
another component is undercrossing a component of the sublink we get 
a piercing of the surface. 

By a {\it minimal spanning surface} for a link we always 
mean an oriented surface in $S^3$, without closed 
simply connected components, 
of maximal Euler characteristic 
that has the link as its boundary. This surface need not
be connected.

\subsection{The Thurston norm}

Let $M$ be a $3$-manifold and $b_i$ the $i$-th Betti number $M$, i.e.\ 
the rank of $H_i(M;\R)$.  In \cite{Thurston:norm} Thurston defined a
semi-norm $\tn{.}$ on $H^1(M; \R)$.  For most interesting cases, e.g.\
for hyperbolic manifolds, this semi-norm is actually a norm.

Given a properly embedded connected surface $S$ in $M$ let 
$$\chi_-(S):=\max \{-\chi(S), 0 \},$$ where $\chi(S)$ is the 
Euler characteristic of
$S$. If $S$ is not connected then $\chi_-(S)$ is defined to be the sum
of $\chi_-(S_i)$ over each connected component $S_i$ of $S$.

Given a class $C$ in $H^1(M;\Z)$ define $\tn {C}:=\min
\{\chi_-(S)\}$ where $S$ is a properly embedded oriented surface which
represents the dual to $C$ in $H_2(M,\partial M;\Z)$.  
The definition
extends uniquely to a semi-norm on $H^1(M;\R)$. This semi-norm is 
symmetric about the origin and the unit-ball is a convex finite-sided
polyhedron with rational vertices. If $\chi_{-}(S)$ equals the
Thurston norm of the homology class $C$ of $S$ then we say that $S$ is
Thurston norm minimizing for $C$. Gabai \cite{Gabai} showed that
one can relax the requirement for the surfaces to be properly embedded
and work with arbitrary images of surfaces instead. 

\begin{figure}[ht!]
\centerline{\small
\SetLabels 
\E(0.33*-0.2){$\mu_i$}\\
\E(1.05*0.5){$L_i$}\\
\endSetLabels 
\AffixLabels{\scalebox{0.7}{\includegraphics{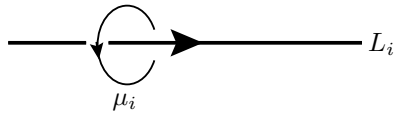}}}}
\vspace{2mm}
\caption{\label{orientation} Orientation of the meridians}
\end{figure}

We are mainly interested in $3$-manifolds that are link complements in
$S^3$.
For a smoothly embedded oriented link $L = L_1 \cup \dots
\cup L_r \subset S^3$ with $r$ components, let $M = S^3 - {\cal N}(L)$ be
the compact $3$-manifold obtained by deleting a tubular neighborhood
${\cal N}(L)$ of $L$.  The orientation of $L$ induces an orientation
on the meridians (see Figure \ref{orientation}) 
and thus, after choosing an ordering of the
components, a multiplicative basis $\langle t_1, \dots, t_r \rangle$
for $H_1(M;\Z) \cong ab(\pi_1(M))$.

Let $\mu_i$ be a meridian curve in $\partial {\cal N}(L_{i})$ 
representing $t_i$.
For a class $C = (C_{1}, \ldots, C_{r}) \in H^{1}(M ; \Z)$, with 
respect to the basis dual to the $t_{i}$, we get:

\begin{lemma}
\label{boundary}
There exists a Thurston norm minimizing surface $S$ for the class $C$ 
such that $S$ is a properly embedded oriented surface in $M$ with $S 
\cap \partial {\cal N}(L_{i})$ a collection $\gamma_{i} = \{ \gamma_{i_1}, 
\dots, \gamma_{i_k} \}$ of parallel simple closed curves, essential in 
$\partial {\cal N}(L_{i})$, such that $|\mu_{i} \cap \gamma_{i}| = 
|C_{i}|$.
\end{lemma}

\begin{proof}
That there exists a properly embedded oriented norm minimizing surface 
$S$ for any class $C \in H^{1}(M; \Z)$ is explained in \cite{Thurston:norm}.  
It remains to show that we may arrange to have 
$\partial S$ as indicated.

We may assume that $S$ is incompressible and boundary incompressible, 
since compressions do not change the homology class of the surface 
and do not raise the value of $\chi_{-}(S)$.  So, since $\partial M$ is a 
collection of tori, $S$ has $S \cap \partial {\cal 
N}(L_{i})$ consisting of some number (possibly zero) of parallel 
simple closed curves, essential in $\partial {\cal N}(L_{i})$.

If there are two components of $\gamma_{i}$ with opposite orientations, 
then there is an innermost pair of such components.  That is, there 
are two components $\gamma_{i_1}$ and $\gamma_{i_2}$ (without loss of 
generality) with opposite orientations such that $\gamma_{i_1} \cup 
\gamma_{i_2}$ is the boundary of an annulus $A$ in $\partial {\cal 
N}(L_{i})$ with $\gamma_{i} \cap \textrm{int}(A)$ empty.  We can 
alter the surface $S$ by attaching a copy of the annulus $A$ to $S$ 
and pushing it into the interior of $M$.  This operation does not 
change the homology class of $S$ or $\chi_{-}(S)$.  By repeatedly 
applying this procedure, we arrange that every component of $S \cap 
\partial {\cal N}(L_{i})$ has the same orientation.

Now homotop $\mu_{i}$ so that $\mu_{i} \cap \gamma_{i}$ is minimal.  Then 
the absolute value of the algebraic intersection number of $\mu_{i}$ 
and $\gamma_{i}$ equals the number of points of intersection of $\mu_{i}$ 
and $\gamma_{i}$.  That is, $|C_{i}| = |\mu_{i} \cap \gamma_{i}|$.
\end{proof}

\subsection{The Alexander norm and McMullen's inequality}

To get an estimate for the Thurston norm, McMullen defined a norm which 
is related to a
classical invariant in knot theory: The Alexander norm $\an{.}$ on
$H^1(M;\Z)$.  For this, one looks at the multivariable Alexander
polynomial for a $3$-manifold \cite{Turaev:Alexanderpolynomial}.  For
$t_1, \dots, t_r$ a set of generators of $H_1(M; \Z)/\mbox{(torsion)}$,
it is a polynomial $\Delta_M$ in $\Z[t_1, \dots, t_n]$. We are only
dealing with link complements, so the $t_i$ can be assumed to be
homology classes of the meridians of the link.

For a class $C \in H^1(M;\Z) \cong Hom(H_1(M;\Z),\Z)$, the Alexander
norm of $C$ is defined to be $\an{C} := \sup C (P - Q)$ where $P$ and
$Q$ ranges over those monomials $t_{j_1} \cdots t_{j_l}$ which appear
with non-zero coefficients in the Alexander polynomial of $M$.  Again,
this norm is sometimes only a semi-norm.

McMullen \cite{McMullen:ThurstonAlexanderNorm} proved an inequality
for the Thurston norm in terms of the Alexander norm.  
For all $C \in H^1(M;\Z)$, $M$ a link complement,
\begin{eqnarray}
\an {C}&\leq& \tn{C}  \qquad \qquad \mbox{for }b_1(M) \geq 2.
\end{eqnarray}
If $M$ is a knot complement and $C$ the generator of $H^1(M)$ we have:
\begin{eqnarray}
\an {C}&\leq& \tn{C} + 1. 
\label{Bettinumber1equation}
\end{eqnarray}
  Inequality
(\ref{Bettinumber1equation}) is the well-known
fact that the span of the Alexander polynomial of a knot is less than or
equal to twice its genus.
If $M$ is not a link complement then one has to take the third
Betti number into account in Inequality (\ref{Bettinumber1equation})
\cite{McMullen:ThurstonAlexanderNorm}.

\subsection {Bennequin's inequality}
\label{sec:Benn_inequality}

There is another, well-known, inequality for the Thurston norm
at certain points in $H^{1}(M; \mathbb R)$. In the sequel of this 
section, we will extend it to all points.
Let $\br \in B_n$ be a braid that closes to a link $L$, where the
orientation is so that all strands are oriented in the same direction.

\begin{figure}[ht]
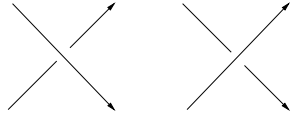

\begin{center} 
\parbox{2 cm}{\scalebox{0.37}{\includegraphics{poscross.eps}}}
\,
\parbox{2cm} {\scalebox{0.37}{\includegraphics{negcross.eps}}}
\end{center}
\caption{\label{signconvention} Positive and negative crossing}
\end{figure}

The {\it Bennequin number} $\beta_t$ for $\br$ is defined to be
$$
\beta_t(\br):=pos(\br) - neg(\br) -n
$$
where $pos(\br)$ is the number of positive crossings and
$neg(\br)$ the number of negative crossings (see Figure \ref {signconvention}).
The difference  $pos(\br) - neg(\br)$ is also called {\it algebraic
crossing number}. (The $t$ in $\beta_t$ stands for ``transversal'' since
it is closely related to and has its origins in transversal knot theory, see 
e.g.\ \cite{FuTa:Legendrian, BirmanWrinkle:transversal}.)

Since, by the algorithm of Pierre Vogel \cite{Vogel}, every link
can be transformed into a closed braid just by Reidemeister II moves
without changing the number of Seifert circles, 
this definition extends to arbitrary link diagrams by replacing $n$ by
the number of Seifert circles in the link diagram.

Let $S$ be an oriented Seifert surface without closed connected 
components that spans the closure of $\br$ and maximizes the 
Euler characteristic $\chi$. 

In his thesis and the paper thereafter \cite {Bennequin} Bennequin
proved:

\begin {theorem} [Bennequin's inequality]
  $\beta_t(\br) \leq -\chi(S)$.
\end{theorem}   

By definition, $-\chi(S) \leq \chi_-(S)$ for a surface $S$.  Hence, 
Bennequin's inequality gives us a lower bound for the Thurston norm at
the points $C=\{C_1, \dots, C_r\}$ with $\vert C_i \vert =1$ for all
$i$.

Note that the Bennequin number depends on the form of the realization of
the link as a closed braid and is not an invariant of the link itself.
The question of finding the maximal Bennequin number among all representatives
of a given link type seems to be harder. We will come back to this
problem in Section \ref{section HOMFLY}.

\subsubsection{Relative Bennequin number} 

For a link, given as a closed braid or as a diagram, 
the Bennequin number gives a lower bound
for the Thurston norm of a single cohomology class. The Bennequin
numbers for sublinks give a much weaker lower bound for
corresponding cohomology classes.
Our aim is to find a suitable definition for a relative Bennequin
number to strengthen these lower bounds and provide
lower bounds for all cohomology classes. 

We proceed as follows:
Given a link $L= L_1 \cup \dots \cup L_r$ and a class $C=(C_1,\dots,C_r)$
in $H^1(M; \Z)$, where $M$ is the link complement, we want to get
a lower bound for the Thurston norm which extends the Bennequin
inequality. 

\begin{itemize}
\item First, we change the orientation of $L$ so that the corresponding
change of the cohomology basis makes all $C_i \geq 0$.
\item Now we transform this reoriented link $L$ into a closed braid,
using Vogel's algorithm. This procedure is by no means unique. 
All our constructions and inequalities from now on will depend on the 
choice of this braid representative for $L$. So we
assume that all $C_i \geq 0$ and in addition the link is given as
a closed braid.
\item In a first theorem we will prove an inequality if all the
$C_i$ are either $1$ or $0$. This corresponds to the question of
finding a minimal spanning surface of the sublink of those components
of $L$ that correspond to $C_i=1$ in the complement of the others.
\item The final step is to explicitly construct a link diagram 
so that the Thur\-ston norm minimizing surface which corresponds to the class
$C$ is a minimal spanning surface for that link. Here we can apply the tools
that we worked out before to get a lower bound for all $C$ with $C_i \geq 0$.
\end{itemize}

We will call the lower bound the {\it relative Bennequin number}.
The definition of the relative Bennequin number of subdiagrams of
link diagrams is meant to be an extension of the Bennequin number 
taking into account the crossings of components not in the sublink.
Let $cr(L_{i},L_{j})$ denote the algebraic crossing number between 
components $L_{i}$ and $L_{j}$ of $L$ and $n_{i}$ the number of 
Seifert circles in Seifert's algorithm applied to $L_{i}$ alone.
We define the {\it relative Bennequin number} $$\beta_t(L_{i},L):= 
cr(L_{i},L_{i}) - n_{i} +
\sum_{j \ne i}\frac{1}{2}cr(L_{i},L_{j}) = cr(L_{i},L_{i}) - n_{i} +
\sum_{j \ne i}lk(L_{i},L_{j})$$ 
and define $$
\beta_t(L_{i_{1}} \cup \dots \cup L_{i_{k}},L) := \sum_{j=1}^{k} 
\beta_t(L_{i_{j}},L).$$  
For a closed braid, note that $\beta_t(L,L) = \beta_{t}(L)$ since the 
number of Seifert circles is the sum of the numbers of Seifert circles 
of the individual components.

We can extend Bennequin's inequality to an inequality between the 
relative Bennequin number and the Thurston norm.

\begin{theorem}
\label{relbennequin}
Let $1 \le i_{1} < i_{2} < \dots < i_{k} \le r$ be integers and let $C 
= (C_{1}, \dots, C_{r})$ where $C_{i}=1$ if $i \in \{i_{1}, \dots 
i_{k} \}$ and $C_{i}=0$ otherwise.  Then
$$\beta_{t}(L_{i_{1}} \cup \dots \cup L_{i_{k}},L) \leq \|C\|_{T}.$$
\end{theorem}

\begin{proof}
Suppose that $S$ is a norm minimizing surface representing $C$.  By 
Lemma \ref{boundary} we may assume that $S \cap \partial {\cal 
N}(L_{i_{j}})$ is a single simple closed curve isotopic to $L_{i_{j}}$ 
for each $1 \le j \le k$ and that $S \cap \partial {\cal N}(L_{i})$ is 
no fewer than $|\sum_{j=1}^{k}lk(L_{i},L_{i_{j}})|$ meridional curves 
for each $i \notin \{i_{1}, \dots i_{k} \}$.
  
We can fill in each boundary curve of $S$ on $\partial {\cal 
N}(L_{i}), i \notin \{i_{1}, \dots i_{k} \}$ with a meridian disk for 
${\cal N}(L_{i})$ to obtain a Seifert surface $S'$ for the link 
$L_{i_{1}} \cup \dots \cup L_{i_{k}}$ in $S^{3}$.  Then $$-\chi(S') 
\leq -\chi(S) - \sum_{i \notin \{i_{1}, \dots, i_{k} 
\}}|\sum_{j=1}^{k}lk(L_{i},L_{i_{j}})|.$$

Bennequin's inequality tells us that $\beta_{t}(L_{i_{1}} 
\cup \dots \cup L_{i_{k}}) \leq -\chi(S')$.  So 
\begin{eqnarray}
\beta_{t}(L_{i_{1}}\! \cup \dots \cup L_{i_{k}},L)\!\!\! &=& 
\beta_{t}(L_{i_{1}}\! \cup \dots \cup L_{i_{k}}) + \! \! \! \sum_{i \notin 
\{i_{1}, \dots, i_{k} \}} \sum_{j=1}^{k} lk(L_{i},L_{i_{j}}) \nonumber 
\\
& \leq & 
\beta_{t}(L_{i_{1}}\! \cup \dots \cup L_{i_{k}}) + \! \! \! \sum_{i \notin 
\{i_{1}, \dots, i_{k} \}} \! |\sum_{j=1}^{k} lk(L_{i},L_{i_{j}})| 
\label{betterbound} \\
& \leq & -\chi(S') + \sum_{i \notin \{i_{1}, \dots, i_{k} \}} 
|\sum_{j=1}^{k}lk(L_{i},L_{i_{j}})| \nonumber \\
& \leq & -\chi(S) \nonumber \\
& \leq & \chi_{-}(S). \nonumber
\end{eqnarray}\vglue -30pt
\end{proof}

\begin{remark}
\label{remarkbetterbound}
In fact, we have proved a somewhat better 
bound holds.  In the proof of Theorem \ref{relbennequin}, we showed that 
the quantity in line (\ref{betterbound}) is a lower bound for the 
Thurston norm of the class $C$.  This inequality is stronger than the 
inequality in the statement of the theorem if the total linking of 
any component $L_{i}$ for $i \notin \{i_{1}, \dots i_{k} \}$ with the 
components $L_{i_{j}}$ is negative.
\end{remark}

We would like to make an assertion like Theorem \ref{relbennequin} for 
any class $C= (C_1, \dots, C_r) \in H^{1}(M;\Z), C_i \geq 0$.  
To do so, we must associate a relative 
Bennequin number to classes with coordinates greater than 1.  We will 
use the closed braid diagram for $L$ to 
create a new closed braid link diagram $L'$ and a subdiagram $L''$ so 
that the pair corresponds to the class $C$.  These new link diagrams 
are derived from the boundary of a Thurston norm minimizing 
surface $S$ for $C$.  More precisely, $L'$ is the link constructed by 
replacing each component $L_{i}$ of $L$ where $C_{i} > 0$ with the 
boundary components of $S$ which lie on the boundary of a regular 
neighborhood of $L_{i}$.  While the complements of $L$ and $L'$ are 
typically not homeomorphic, the surface $S$ exists in both.  By 
passing from $L$ to $L'$, we will have arranged for $S$ to be dual to an 
element in $H^{1}(S^{3} - {\cal N}(L'))$ that is of the form where we can 
apply the classical Bennequin inequality.

By Lemma \ref{boundary} we may assume that a norm minimizing surface 
$S$ representing $C$ has $S \cap \partial {\cal N}(L_{i})$ a 
collection $\gamma_{i}$ of parallel simple closed curves that are 
essential in $\partial {\cal N}(L_{i})$ and all oriented in the same 
direction.  Let $\lambda_{i}$ be a longitude for $L_{i}$, that is, a 
simple closed curve in $\partial {\cal N}(L_{i})$ such that 
$\lambda_{i}$ does not 
link $L_{i}$ and such that the intersection number of $\mu_{i}$ 
and $\lambda_{i}$ is equal to $+1$.
We have that $[\gamma_{i}] = p_{i}[\mu_{i}] + q_{i}[\lambda_{i}] 
\in H_{1}(\partial {\cal N}(L_{i}) ; \Z)$ where 
$$q_{i} = C_{i} \ \mbox{and} \  p_{i} = -\sum_{j 
\ne i} C_{j} \cdot lk(L_{i},L_{j}).$$
We define the two new link diagrams based on the closed braid diagram 
for $L$.  Let $L'$ be the link diagram obtained by replacing the link 
component $L_{i}$ with the curves $\gamma_{i}$ for all $1 \le i \le r$ 
with $C_{i} \ne 0$ as follows.  Replace each strand of the component 
$L_{i}$ by $C_{i}$ parallel strands.  Let the new strands replacing 
the first strand of $L_{i}$ be strands $m_{i}+1, m_{i}+2, \dots, 
m_{i}+C_{i}$.  At the end of these strands, insert the braid 
$$(\sigma_{m_{i}+1} \sigma_{m_{i}+2} \dots 
\sigma_{m_{i}+C_{i}-1})^{p_{i} - C_{i} cr(L_{i},L_{i})}.$$
Here, the $\s_i$ denote a positive standard generator of the braid
group, i.e.\ the braid has a positive crossing between the $i$-th. and
$i+1$-th. strand.
We obtain 
the link diagram $L''$ from $L'$ by deleting each component $L_{i}$ 
with $C_{i}=0$.  Note that $L'$ and $L''$ are closed braids and that 
$L''$ is a subdiagram of $L'$.

\begin{example}
Let $L$ be the closure of the braid $\s_1^4$, which is
a $2$-component link $L_1 \cup L_2$. Seifert's algorithm gives us a 
torus with two boundary components. Hence, the Thurston
norm at the point $(1,1)$ is less than or equal to $2$.

Indeed, Bennequin's inequality computes as
$$ \beta_t(L)= \beta_t(L_1 \cup L_2) = 4 -2 =2 \leq -\chi(S(L)).$$
Hence, the Thurston norm is $2$ at this point.

The relative Bennequin numbers $\beta_t(L_i,L), i=1,2,$ are
both $1$. Again, since the generalized Seifert algorithm gives
us a two-punctured disk as a spanning surface of one component in
the complement of the other, we get the Thurston norm at $(1,0)$ and
$(0,1)$.

Now, look for example at the point $(2,1)$. 
In the construction above we get $q_1 = 2, q_2 = 1$ and $p_1 = -2$.
Hence, our new closed braid  
$L' = L''$ is $\s_2 \s_1^2 \s_2^2 \s_1^2 \s_2 \s_1^{-2}$
and the relative Bennequin number $\beta_t(L',L')$ equals $3$.

The Thurston norm is a semi-norm and we get
$$3 = \beta_t(L',L') \leq \tn{(2,1)} \leq \tn{(1,1)} + \tn{(1,0)} = 2 +1 =3.$$
\end{example}

\begin{theorem}
\label{additive}
$$\beta_{t}(L'',L') = \sum_{i=1}^{r} C_{i} \beta_{t}(L_{i},L).$$
\end{theorem}

\begin{proof}
We first note that $\beta_{t}(L'',L') = \beta_{t}(L'') + 1/2 cr(L'', 
L' - L'')$ where $cr(L'', L' - L'')$ is the algebraic crossing number
between components in $L''$ and components not in $L''$.

Crossings in $L''$ occur in the following ways:
\begin{itemize}

\item For each crossing of $L_{i}$ with itself, there are $C_{i}^{2}$ 
crossings with the same sign in $L''$.

\item For each component $L_{i}$, there are $(C_{i}-1)(p_{i} - C_{i} 
cr(L_{i},L_{i}))$ crossings in $L''$ (counted with sign) coming from 
the braid inserted above.

\item For each crossing of $L_{i}$ with $L_{j}$ for $j \ne i$, there 
are $C_{i} C_{j}$ crossings with the same sign in $L''$.

\end{itemize}

Therefore, 
\begin{eqnarray*}
\beta_{t}(L'') & = & \left . \sum_{i=1, C_i \neq 0}^{r} 
\right ( C_{i}^{2} cr(L_{i},L_{i}) - 
C_{i}(C_{i}-1) cr(L_{i},L_{i}) \\
& & \left . + p_{i} (C_{i}-1) + \sum_{j \ne 
i} \left (C_{i} C_{j} lk(L_{i}, L_{j}) \right ) - C_{i} n_{i} \right ) \\
& = & \left . \sum_{i=1, C_i \neq 0}^{r} \right ( C_{i}^{2} cr(L_{i},L_{i}) - 
C_{i}(C_{i}-1) cr(L_{i},L_{i} )  \\
& & - (C_{i}-1) \left . \sum_{j \ne i} (C_{j} lk(L_{i},L_{j})) + \sum_{j \ne 
i} (C_{i} C_{j} lk(L_{i}, L_{j})) - C_{i} n_{i} \right ) \\
& = & \sum_{i=1, C_i \neq 0}^{r} \left (C_{i} cr(L_{i},L_{i}) + \sum_{j \ne i} C_{j} 
lk(L_{i},L_{j}) - C_{i} n_{i}\right ).
\end{eqnarray*}
The sums can be taken over all $j$ in these equations because 
each relevant term contains a $C_{j}$ factor, 
eliminating those components in $L' - L''$.

Hence,
\begin{eqnarray}
\lefteqn{\beta_{t}(L''\!,\! L')=} \nonumber\\
 &=& \beta_t(L'') + 1/2 cr(L'',L'-L'') \nonumber \\ 
& = &\!\! \sum_{i=1,C_i \neq 0}^{r} \left(C_{i} cr(L_{i},L_{i}) + \sum_{j 
\ne i} C_{j} lk(L_{i},L_{j}) - C_{i} n_{i}\right) + 1/2 cr(L'', L' - L'') \nonumber \\
& = &\!\! \sum_{i=1, C_i \neq 0}^{r} \left(C_{i} cr(L_{i},L_{i}) + \sum_{j \ne i} C_{j} 
lk(L_{i},L_{j}) - C_{i} n_{i} + 1/2 C_{i} \sum_{C_{j}=0}cr(L_{i}, 
L_{j})\right) \nonumber \\
& = &\!\! \sum_{i=1, C_i \neq 0}^{r} \left(C_{i} cr(L_{i},L_{i}) + \sum_{j \ne i} C_{j} 
lk(L_{i},L_{j}) - C_{i} n_{i} + C_{i} \sum_{C_{j}=0}lk(L_{i}, 
L_{j})\right) 
\label{longline} \\
& = & \sum_{i=1}^{r} \left(C_{i} cr(L_{i},L_{i}) + C_{i }\sum_{j \ne i} 
lk(L_{i},L_{j}) - C_{i} n_{i}\right) \label{hardtosee} \\
& = & \sum_{i=1}^{r} C_{i} \beta_{t}(L_{i},L). \nonumber 
\end{eqnarray}
The equality in line (\ref{hardtosee}) is true because the second term 
in line (\ref{longline}) will produce all terms of the form 
$C_{i}lk(L_{i},L_{j})=C_i lk(L_j,L_i)$ for all pairs $(i,j)$ 
with $i \ne j$ and neither $C_{i}$ nor 
$C_{j}$ equal to zero.  The fourth term in line (\ref{longline}) will 
produce all terms of the form $C_{i}lk(L_{i},L_{j})$ for all pairs 
$(i,j)$ with $C_{i} \ne 0$ and $C_{j} = 0$.  Together, they produce 
all terms of the form $C_{i}lk(L_{i},L_{j})$ for all pairs 
$(i,j)$ with $i \ne j$ and $C_{i} \ne 0$.  The sum in 
line (\ref{hardtosee}) can be taken over all $i$ because a factor of 
$C_{i}$ appears in every term, eliminating those with $C_{i}=0$.
\end{proof}

With this linearity result on relative Bennequin numbers, we obtain 
a lower bound on the Thurston norm of any cohomology class $C=(C_1, 
\dots, C_r)$ in terms of the relative Bennequin numbers of each link 
component of the closed braid diagram for $L$ reoriented so that all 
$C_i \geq 0$.

\begin{corollary}
\label{bennequin_ineq}
For any class $C=(C_1, \dots, C_r) \in H^1(M;\Z)$,
$$\beta_t(L'',L') = \sum_{i=1}^{r}C_i \beta_{t}(L_{i},L) \leq \|C\|_T.$$
\end{corollary}

\begin{proof}
The pair of link diagrams $L'$ and $L''$ are defined in terms of 
a norm minimizing surface $S$ for the class $C$ so that $S$ is also 
a surface spanning $L''$ in the complement of $L'$.  Therefore, by 
Theorem \ref{relbennequin}, $\beta_t(L'',L') \leq \chi_-(S) = \|C\|_T$.
\end{proof}

As in Remark \ref{remarkbetterbound}, our proof of Theorem 
\ref{relbennequin} provides a better lower bound than the one given 
in Corollary \ref{bennequin_ineq}, although it does not have as 
succinct an expression.  This lower bound is outlined in the Scholium 
below.

\begin{scholium}
For any class $C=(C_1, \dots, C_r) \in H^1(M;\Z)$,
$$\beta_t(L'',L') \leq \beta_{t}(L'') + \sum_{i=1,C_{i} = 0}^{r} | 
\sum_{j=1}^{r} C_{j} lk(L_{i},L_{j}) | \leq \|C\|_T.$$
\end{scholium}

\proof
The first inequality is actually an equality if each term in the 
absolute value signs is non-negative.  In this case the statement 
here is identical to Corollary \ref{bennequin_ineq}.

As in the proof of Theorem \ref{relbennequin}, a Thurston norm 
minimizing surface $S$ with no closed components representing $C$ has 
$S \cap \partial {\cal N}(L_{i})$ equal to no fewer than $| 
\sum_{j=1}^{r} C_{j} lk(L_{i},L_{j}) |$ meridian curves for each $i 
\in \{1, \dots, r\}$ with $C_{i}=0$.  Filling in each of these with a 
meridian disk produces a spanning surface $S''$ for $L''$ in $S^{3}$.  
Then $$-\chi(S'') \leq -\chi(S) - \sum_{i=1,C_{i} = 0}^{r} | 
\sum_{j=1}^{r} C_{j} lk(L_{i},L_{j}) |.$$ From Bennequin's inequality, 
we know that $\beta_{t}(L'') \leq -\chi(S'')$.  Therefore, $$\beta_{t}(L'') 
+ \sum_{i=1,C_{i} = 0}^{r} | \sum_{j=1}^{r} C_{j} lk(L_{i},L_{j}) | 
\leq - \chi(S) \leq \chi_{-}(S) = \|C\|_T.\eqno{\qed}$$

\section[Closed homogeneous braids] 
{The Comparison of the Thurston and the Alexander norm}

\subsection{Closed homogeneous braids}

First we need a lemma for arbitrary semi-norms

\begin{lemma} \label{additivity}
Let $\| . \|$ on $\R^n$ with basis $e_i, i=1,\dots,n$,
be a semi-norm, such that 
$$
\left \| \sum_{i=1}^n e_i \right \| = \sum_{i=1}^n \| e_i \|.
$$
Then
$$\left \| \sum_{i=1}^n c_i e_i \right \| = \sum_{i=1}^n c_i \| e_i \|$$
for all $c_i \geq 0$.
\end{lemma}
\begin{proof}
Assume that the basis $e_i$ is chosen so that $\| e_i \| \neq 0$ for
$i \leq n_1$ and $\| e_i \| =0$ for $n_1 < i \leq n$.
The unit ball $B$ for $\| . \|$ is convex. Therefore
$B$ is an intersection of half-spaces. If a point $v \in \R^n$ has
norm $\| v \| =1$ then there exists an $(n-1)$-dimensional affine
space ${\cal A}$ through $v$ such that $B$ is contained in the closure of 
one of two halfspaces defined by ${\cal A}$. Let $v_i$ be the point 
$e_i/\| e_i\|$ for $i \leq n_1$. Choose $v = (\sum e_i)/(\sum \| e_i 
\|)$, hence $\| v \| =1$. The only $(n-1)$-dimensional
affine space ${\cal A}$ that contains $v$, and has a closed halfspace
that contains all  $v_i, i \leq n_i,$ and all $m_i e_i, i >n_1,$ is
the $(n_1-1)$-dimensional affine space  
containing the $v_i$'s times $R^{n-n_1}$. 

For all  $w \in {\cal A} \cap \R_{\geq 0}^n$ it follows $\| w \| =1$.  
\end{proof}

We get:

\begin{theorem} \label{Thurston norm on positives}
Let $\br$ be a braid so that Seifert's algorithm gives a minimal
spanning surface for its closure $\hat \br$ and let 
$M=S^3 - {\cal N}(\hat \br)$ be the complement. Assume furthermore
that
$\hat \br$ has no unlinked unknotted components. 
Number the components of $\hat \br$ and take as a basis for 
$H^1(M;\Z)$ the one induced by the meridians
as described in Section \ref{setting_the_scene}.

The Thurston norm is additive 
for all $(C_1, \dots, C_r) \in H^1(M;\Z)$ with $C_j \geq 0$.
Its value on the points $(0,\dots,0,C_j=1,0,\dots,0)$ equals minus
the Euler characteristic of the spanning surface for the $j$-th component
pierced by the other components that we get as described in Section
\ref{generalized Seifert algorithm}.
\end{theorem}

\begin{proof}
First we prove the theorem for the case that all $C_i$ are equal to 
either $0$ or $1$.  By assumption the Thurston norm at the point 
$(1,\dots,1)$ equals
$$\tn{(1, \dots,1 )} = l - n$$
where $l$ is the word length of the braid and $n$ the 
number of strands. 

Choose a component $L_{j}$ with $n_j$ strands of the link. Construct
a surface that spans this component as in the generalization of
Seifert's algorithm explained in Section \ref{generalized Seifert algorithm}. 
That is, if another component undercrosses an arc of $L_{j}$, we get a
piercing of the surface. Let $u_j$ be the number of these piercings and
$l_j$ be the number of self-crossings of the component $L_{j}$.
The so constructed surface has Euler characteristic
$$n_j - u_j - l_j,$$ and thus, the Thurston norm at the point
$(0, \dots, 0, C_j=1, 0, \dots 0)$ is less than or equal to $u_j+l_j-n_j$.
By definition $n= \sum n_j$
and $l = \sum (l_j + u_j)$.

Since the Thurston norm is a semi-norm, we have
\begin{eqnarray*}
l - n = \tn{(1,\dots,1)} &\leq& \sum_j \tn{(0,\dots,0,C_j=1,
0,\dots,0)} \\
& \leq&  \sum_j (u_j + l_j - n_j) = l - n.
\end{eqnarray*}
Thus, the Thurston norm of 
$(0,\dots,0,C_j=1,0,\dots,0)$ equals $l_j+u_j-n_j$.

Lemma \ref{additivity} gives us
the additivity for all other points $(C_1, \dots, C_r)$
with $C_i \geq 0$. 
\end{proof}
 
\begin{proposition} \label{Alexander norm equals Thurston norm}
Assume the conditions for the Thurston norm as in Theorem
\ref{Thurston norm on positives} are fulfilled. Assume
furthermore, that the Alexander norm at the point 
$(1,\dots, 1) \in H^1(M;\Z)$
equals the Thurston norm.

Then the Alexander norm equals the Thurston norm for all
points $(C_1, \dots, C_r)$ with $C_i \geq 0$.
\end{proposition}

\proof
If the closed braid is a knot, then the theorem is immediate.  Hence, 
we assume that $\hat \br$ has at least two components.  Since the 
Alexander norm is bounded from above by the Thurston norm, it follows 
that for $C^{(j)}= (0, \dots, 0, C_j=1, 0, \dots, 0)$, $$\an{(0, 
\dots, 0, C_j=1, 0, \dots, 0)}=l_j+u_j-n_j = \tn{(0,
\dots, 0, C_j=1, 0, \dots, 0)}.$$
Now, Lemma \ref{additivity} shows that the Alexander norm is also 
additive on all points $C = (C_1, \dots, C_r)$ with $C_i \geq 0$. 
Therefore, $$\an{C} = \sum_{j=1}^{r} C_{j} \an{C^{(j)}} = \sum_{j=1}^{r} 
C_{j} \tn{C^{(j)}} = \tn{C}.\eqno{\qed}$$
An important class of links are closures of homogeneous braids,
i.e.\ the braid contains every standard generator 
$\s_i$ of the braid group at least once and always with the 
same exponent sign.

\begin{corollary}
For closed homogeneous braids the Alexander norm equals the Thurston
norm for all points $(C_1, \dots, C_r) \in H^1(M;\Z)$ with $C_i \geq 0$.
\end{corollary}

\begin{proof}
By a theorem of Stallings \cite{Stallings:fibred_links} these
links are fibred and the fibre surface is the one that one gets
by the application of Seifert's algorithm to the closed braid.

Hence, the Alexander norm equals the Thurston norm at the point
$(1, \dots, 1)$ (e.g.\ \cite{McMullen:ThurstonAlexanderNorm}) 
and the claim follows by Proposition 
\ref{Alexander norm equals Thurston norm}.
\end {proof}

\section {The Bennequin number and the Thurston norm}

\subsection{The Bennequin number and the band generators}
\label{band generators}

\begin{figure}[ht!]
\centerline{\small
\SetLabels 
\E(-0.05*0.15){$i$}\\
\E(-0.05*0.85){$j$}\\
\endSetLabels 
\AffixLabels{\scalebox{0.6}{\includegraphics{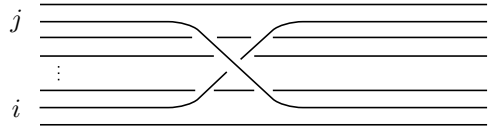}}}}
\caption{\label{fig band generators} The band generator $a_{i,j}$}
\end{figure}

Recently, a new presentation for the braid groups began to receive attention:
a presentation in terms of the band generators, defined by Birman, Ko and Lee 
\cite{BKL}. Let 
$$a_{i,j}:= \s_i^{-1} \cdots \s_{j-2}^{-1} \s_{j-1} \s_{j-2} \cdots
\s_i,$$
where $\s_i$ are the standard generators 
(see Figure \ref{fig band generators}).
In \cite{BKL} a solution to the word problem and conjugacy problem 
with respect to these band generators is given.

For a closed braid given in terms of the band generators,
we get a spanning surface by a generalized Seifert algorithm:
the number of Seifert disks equals the number of strands, and
if $a_{i,j}$ or $a_{i,j}^{-1}$ occurs in a word, then we connect the
$i$-th and $j$-th disk by a band.
Thus, the Euler characteristic of this spanning surface is $n -l$,
where $l$ is the word length in the band generators.

For $\br \in B_n$ a braid, let
$\pos_b(\br)$ be the positive exponent sum in the band
generators and $\negat_b(\br)$ be the negative exponent sum.
Bennequin's inequality gives us $$\pos_b(\br)-\negat_b(\br) -n \leq 
-\chi(S(\hat \br))$$
for any spanning surface of the closed braid $\hat \br$.

A braid $\br$ is band positive if $\negat_b(\br)=0$.
(This is a different terminology than in the work of Rudolph
(see e.g.\ \cite{Rudolph} and the references therein).
There, band positive braids are a special type of
quasi-positive braids.)

Bennequin's inequality gives the following proposition which is an 
extension of the main
result of Ko and Lee in \cite{KL:genera}, where it is proved, with
different methods, if 
the number of strands equals $4$.

\begin{proposition} \label{band positive}
Let $\br$ be a band positive braid in $B_n$ of length $l$.
Then the Euler characteristic of a minimal spanning surface 
satisfies:  $-\chi(S(\hat \br))= l-n$.
\end{proposition}

A special class of band positive braids is particularly well-suited
for computing the Thurston norm at the cohomology classes 
$C=(C_1,\dots,C_r)$ with $C_i \geq 0$.  This class consists of band
positive braids with only bands of the 
form $a_{i,i+1}$ or $a_{i,j}$ 
such that the $i$th and $j$th strands are parts of the same link component
in the closed braid.  
This class of braids has the property that every sublink is also band 
positive.

\begin{proposition}
For a closed braid in the class described above with no unlinked, unknotted 
components, 
the relative Bennequin inequality 
in Corollary \ref{bennequin_ineq} is an equality for every cohomology
class $C=(C_1,\dots,C_r)$ with $C_i \geq 0$.
\end{proposition}

\begin{proof}
First, on the points $C^{(i)}=(0, \dots, 0, C_i =1, 0, \dots,0)$ the
relative Bennequin inequality is $\beta_t(L_i,L) \leq \tn{C^{(i)}}$.
On the other hand, a generalized Seifert band algorithm similar to
the generalized Seifert algorithm for the standard generators described
in Section \ref{generalized Seifert algorithm} produces a surface $S_i$ with 
$\chi_{-}(S_i)= -\chi(S_i) = \beta_t(L_i,L)$.

Hence, we know that the relative Bennequin inequality is an equality for
all points $C^{(i)}$ and in addition for
$(1,\dots,1)$. Thus, by Lemma \ref{additivity}, the Thurston norm is 
additive on
all points $(C_1, \dots, C_r)$ with $C_i \geq 0$.
Since we know by Theorem \ref{additive} that the relative Bennequin 
number
is also additive on these points, we are done.  
\end{proof}

Note that, in contrast, the estimate on the Alexander norm coming from the 
Thurston norm for this class of closed braids is in general not sharp.
In \cite{KL:genera} an example of a knot represented by a 
closed band positive braid is given, where the Thurston
norm equals $3$, while the Alexander norm is $2$.

Using this example one can easily find more examples of links with 
more components. The computation was made with the help of the
program  \cite{Morton:Computerprogram_MVA} based on a description
of the multivariable Alexander polynomial due to Morton
\cite {Morton:MVA}.

The two component link represented by the closure of the 
$5$-braid $$a_{4,5}^2 a_{2,4}^2 a_{1,3} a_{3,4} a_{2,4} a_{1,3}^2$$ 
has Alexander polynomial $2 -3 t_1 + 2 t_1^2$ (the second variable
vanishes), hence the Alexander norm at $(1,1)$ is $2$.
The Thurston norm equals $4$, by Proposition \ref{band positive}. 

\section{The HOMFLY polynomial and a conjecture of\break Morton}
\label{section HOMFLY}

In the last section we deal with another invariant for links which
is conjectured to also give a lower bound for the Thurston norm.
It comes from the HOMFLY polynomial.
Shortly after the HOMFLY polynomial was defined, Hugh Morton 
\cite{Morton:Seifert_circles} and independently
Franks and Williams \cite{FW:Bennequin_estimate} 
realized that its lowest degree in one of the variables
can be estimated by the Bennequin number (see e.g.\ 
\cite{Lickorish:knotpolynomials} for a still excellent account to the, by now,
classical knot polynomials).

We briefly recall the relevant results. On braids let 
a two-variable (Laurent-) polynomial $P(\br)(v,z)$ be defined by:

\begin{eqnarray}
P(\br \s_i)-P(\br \s_i^{-1})&=& z P(\br), \qquad 1 \leq i \leq n-1\label{skein}\\ 
P(\br \s_n) &=& P(\br) \label{pos_stab}\\
P(\br \s_n^{-1}) &=& v^2 P(\br) \label{neg_stab}\\
P(\id_n)&=& ((1-v^2)/z)^n \label{loop}  
\end{eqnarray}
Here $\br$ is  assumed to be a braid in $B_n$, thus $\s_n$ does not
occur in $\br$, and
$\id_n$ is the identity in $B_n$.
Furthermore, the polynomial is, by definition, 
invariant under braid isotopy, i.e.\ relations
in $B_n$, and conjugation.
From the definition it follows that if the braid has the form
$\id_n \br$ with $\br \in B_{n-1}$ then
\begin{eqnarray} \label{equation circle in HOMFLY}
P(\id_n \br) &=& ((1-v^2)/z) P(\br).
\end{eqnarray}
It is easy to see, that 
$$v^{\BeT(\br)} \, P(\br)=\frac{(v^{-1}-v)}{z} \HOMFLY(\hat \br)$$
where $\HOMFLY(\hat \br)$ is the HOMFLY polynomial of the closure $\hat \br$ of
$\br$ and $\BeT(\br)$ is again the Bennequin number of the braid $\br$,
i.e.\ $\pos(\br)-\negat(\br)-n$.
This gives that $P$, assuming the definitions are sufficient,
is unique and well-defined.
To see that the definition is sufficient to compute $P$, notice that
\begin{eqnarray}
P(\br \s_i^2)&=& P(\br) + z P(\br \s_i) \label{positive reduction}  \\
P(\br \s_i^{-1}) &=& P(\br \s_i) - z P(\br). \label{negative reduction}
\end{eqnarray}
With (\ref {negative reduction}) one can assume that it is sufficient to compute 
$P$ for all
positive braids. Relations (\ref{positive reduction}), (\ref{pos_stab}), (\ref{neg_stab})
and (\ref{loop}) are now sufficient to compute $P$ without having to increase 
the word length.

The specialization $v:=1$ in $$z \, P(v,z)(\br)/ (1-v^2)$$ is the 
(one variable) Alexander 
polynomial of the closed braid $\hat \br$ in its Conway form.

From the definition of $P$ and its computation it follows immediately
that $P$ is a polynomial in $v$ rather than a Laurent polynomial.
Therefore, if $e(\hat \br)$ is the minimum degree of $v$ in 
$\HOMFLY(\hat \br)$ we have

\begin{theorem}[Morton \cite{Morton:Seifert_circles} , 
Franks-Williams \cite{FW:Bennequin_estimate}] 
\label{Bennequin-Homfly inequality}
Let $\hat \br$ be the\break closure of a braid $\br$ and $e(\hat \br)$ be the 
minimum degree of $v$ in the HOMFLY polynomial of $\hat \br$. Then
\begin{eqnarray*}
\BeT(\br) +1 \leq e(\hat \br).
\end{eqnarray*}
\end{theorem}

This theorem turned out to be a powerful tool in Legendrian and transversal
knot theory, as an upper bound for the Bennequin number (see
\cite{FuTa:Legendrian,Tabachnikov:Legendrian,CG:Legendrian}).

As a corollary we get:

\begin{corollary} \label{maximal Bennequin number}
If, for a braid $\br$, the polynomial $P(\br)(0,z) \neq 0$ then
the Bennequin number of $\br$ is maximal among the Bennequin numbers of 
all representatives of the 
closure $\hat \br$ of $\br$. 
\end{corollary}

In \cite{Morton:Seifert_circles} Hugh Morton made the intriguing conjecture:
\begin{conjecture}[Morton]
For a link $L$ let $e(L)$ be the lowest degree of the HOMFLY polynomial in
the framing variable $v$. Then
$$e(L) \leq -\chi(S(L)) +1$$ for every Seifert surface $S(L)$ of $L$.
\end{conjecture}

This conjecture and Theorem \ref{Bennequin-Homfly inequality} again would imply the 
Bennequin inequality. 

\subsection{Homogeneous braids}

We discovered that the corollary motivating the 
following proposition can be seen immediately from results 
proven by Cromwell in \cite{Cromwell:homogeneous}. However,
since our approach is braid theoretical rather than based on knot 
diagrams, we include a proof.

\begin{proposition} \label{nontrivial term in P}
Let $\br$ be a homogeneous braid in $B_n$ such that every generator occurs at
least once. Let $n_n$ be the number of
generators $\s_i$ that occur with negative exponents and 
$n_p$ be the number of generators that occur with positive exponents.
Thus, $n=n_n+n_p+1$.
Then the polynomial $P(\br)(v,z)$ contains a unique highest term in $z$:
$$ (-1)^{\negat(\br)-n_n} (1-v^2) v^{2 n_n} \,  z^{\vert \br \vert-n}$$
where $\vert \br \vert$ is the word length of $\br$. 
\end{proposition}

\begin{proof}
We will apply the following elementary equations:
\begin{eqnarray}
P(\br \s_i^2) \!\!&=&\!\! P(\br) + z P(\br \s_i) \label{H transform I}\\
P(\br \s_i^{-2})\!\!&=&\!\! P(\br) - z P(\br \s_i^{-1}) \label{H transform II}\\
P(\br \s_n^{-1} \s_{n-1} \s_n^{-1}) \!\!&=&\!\! P(\br \s_{n-1} \s_n^{-1} 
\s_{n-1}) - z P(\br \s_n^{-1} \s_{n-1}) \nonumber\\
\!\!&&\!\! -z P(\br \s_{n-1} \s_n^{-1})
\label{H transform III}\\
P(\br \s_n \s_{n-1}^{-1} \s_n)\!\!&=&\!\! P(\br \s_{n-1}^{-1} \s_n \s_{n-1}^{-1})
+z P(\br \s_{n-1}^{-1} \s_n)+ z P(\br \s_n \s_{n-1}^{-1}). 
\label{H transform IV} 
\end{eqnarray}
If the proposition is proven for the terms on the right hand side, then it
is also true for the ones on the left hand side.

First, the proposition is true in $B_2$, i.e.\ for braids $\s_1^k$ or $\s_1^{-k}$ by induction on $k$.
For the representative of the unknot $\s_{n-1}^{\epsilon_{n-1}}
\cdots \s_1^{\epsilon_1}, \epsilon_{i} =\pm 1,$ it is also true.

It is enough to prove the claim, that every homogeneous braid in $B_n$ 
can be brought with the transformations 
(\ref{H transform I})-(\ref{H transform IV}) and the braid
relations $\s_{i+1} \s_i \s_{i+1}= \s_i \s_{i+1} \s_i$ 
into one in which
$\s_{n-1}$ occurs at most once, modulo words of lesser length.

For this let the claim already be proven for $B_{n-1}$.
Thus, we can assume that a braid 
$\br:=w_1 \s_{n-1}^{\epsilon_{n-1}} w_2 \s_{n-1}^{\epsilon_{n-1}}$, $w_2 \in B_{n-1}$, modulo terms of lower word length, has in $w_2$ at most one term 
$\s_{n-2}^{\epsilon_{n-2}}$. If $w_2$ is empty then with (\ref{H transform I})
 or (\ref{H transform II}) we can reduce the word length.

Otherwise, since $\s_i$ commutes with $\s_{n-1}$ for $i<n-2$, we can apply
(\ref{H transform III}),(\ref{H transform IV}) or the braid relations 
to bring $\br$ into a form
in which $\s_{n-1}$ occurs less and which is still homogeneous.
\end{proof}

\begin{corollary} {\rm\cite{Cromwell:homogeneous}}\qua
Morton's conjecture is true for closed homogeneous braids.
\end{corollary}

\begin{proof}
By Proposition \ref{nontrivial term in P} and with the notations there, 
we know that the minimum degree
$e(\hat \br)$ of the HOMFLY polynomial is at most $\BeT(\br) + 1 + 2 n_n$.

Since 
\begin{eqnarray*}
\BeT(\br) &=& \pos(\br)-\negat(\br)-n_n-n_p-1  \qquad \mbox{and}\\ 
-\chi(\hat \br) &=& \pos(\br)+\negat(\br) - n_n - n_p -1
\end{eqnarray*}
we have
\begin{eqnarray*}
e(\hat \br) & \leq & \BeT(\br) +1 + 2 n_n\\
& =& -\chi(\hat \br) -2 \negat(\br) + 1 + 2 n_n\\
& \leq & -\chi(\hat \br) + 1.
\end{eqnarray*}
The last equation follows, since $n_n \leq \negat(\br)$. 
\end{proof}

\subsection {Closed $3$-braids}

The aim of this section is to prove Morton's conjecture for
closed $3$-braids. The idea that lies behind the proof is quite 
general, though:

\begin{enumerate}
\item We assume the link to be given as a closed braid $\br$ on $n$ strands
in the band generators, such that the band Seifert algorithm 
(see Section \ref{band generators}) yields a minimal spanning surface $S$.

Thus, we get for the Euler characteristic:
\begin{eqnarray} \label{General:Euler}
-\chi(\hat{\br})=-\chi(S)\!\!&=&\!\! \text{wordlength}(\br)-n=
\pos(\br)+\negat(\br)-n,
\end{eqnarray}
where, again, $\pos(\br)$ ($\negat(\br)$, respectively) is the
number of positive (negative) band generators in the word $\br$.

\item
The Bennequin number computes as
\begin{eqnarray} \label{General:Bennequin}
\BeT(\br)&=& \pos(\br)-\negat(\br)-n.
\end{eqnarray}
\item Morton's conjecture now states that
\begin{eqnarray} \label{General:Morton}
e_P(\br)+\BeT(\br) &\leq& -\chi(\hat{\br}),
\end{eqnarray}
where $e_P(\br)$ is the lowest degree of the polynomial
$P(\br)(v,z)$ in the variable $v$.
\end{enumerate}
Thus, using Equations (\ref{General:Euler}) and (\ref{General:Bennequin})
we can reformulate Equation (\ref{General:Morton}) to get:
\begin {eqnarray*}
e_P(\br)+\pos(\br)-\negat(\br)-n &\leq&\pos(\br)+\negat(\br)-n\\
\iff e_P(\br)&\leq&2 \negat(\br).
\end{eqnarray*}

Since $e_P(\br)\leq 2(n-1)$ by the definition of $P(\br)(v,z)$, it
remains to check Morton's conjecture for all words with
$$\negat(\br) < n-1.$$
In the special case $n=3$ we can prove with this idea:

\begin{theorem} \label{Morton's conjecture for 3-braids}
Morton's conjecture is true for any closed $3$-braid.
\end{theorem}

\begin{proof}
Let again $a_{1,2}:=\s_1, a_{2,3}:=\s_2$ and $a_{1,3}:=\s_1^{-1} \s_2 \s_1$.
With these band generators, the braid group $B_3$ has a presentation
\begin{eqnarray} \label{presentation B3}
B_3 &\cong & \langle a_{1,2},a_{2,3},a_{1,3} \vert a_{1,2} a_{1,3}= a_{2,3} a_{1,2} =
a_{1,3} a_{2,3} \rangle.
\end{eqnarray}
It is known by work of Birman and Menasco \cite {BirmanMenasco:3braids}
and of Xu \cite{Xu} that a Seifert surface with maximal Euler characteristic
is realized by a shortest word in these three generators.
This means, if a closed $3$-braid $\hat \br$ is represented by a shortest 
word of length $l$ 
then a minimal spanning surface $S(\hat \br)$ has $-\chi(S(\hat \br))= l -3$.
It is the surface that we get by the band Seifert algorithm described
in Section \ref{band generators}.

Now Morton's conjecture claims that the Bennequin number of $\br$ 
and $-\chi$ of a minimal spanning surface for $\hat \br$
differ at least by the minimum degree of $v$ in $P(\br)(v,z)$.

By the computation of $P$ it is clear that for a $3$-braid the only 
possible powers of $v$ in the monomials of $P$ are 
$0$, $2$, $4$, and $6$.  Furthermore, $P$ is 
divisible by $((1-v^2)/z)$ and thus the minimum degree of $v$ in 
$P(v,z)$ can be at most $4$.

Now, let $\br$ be a $3$-braid, which is of minimal length, 
with exponent sum $\negat(\br)$ over all negative exponents and
$\pos(\br)$ over all positive exponents.
Thus,  the difference between the Bennequin number and $-\chi$ is 
$2 \, \negat(\br)$. So, we only have to check Morton's conjecture for braids
$\br$ with $\negat(\br) \leq 1$.

First assume, that $\br$ is a positive word in the generators.
Our claim is that the highest degree in $z$ is $l-3$ where $l$ is the
word length. Moreover, we claim that coefficient of $z^{l-3}$ which is
a polynomial in $v$ has constant term $1$.

If the word length is one then, as one might easily check, 
$P(\br)=((1-v^2)/z)^2$.
If the word length is two, the verification of the claim is an easy 
case-by-case check as well.

Assume that the word length is at least $3$. If a square of a generator occurs,
then for any of the generators $a=a_{i,j}$ we have
\begin{eqnarray} \label{induction 1}
P(\br a^2)&=& P(\br) + z P(\br a)
\end{eqnarray}
and the claim holds by induction.
Since we have the relations (\ref{presentation B3}) we can thus assume that
any subword of length three is of the form $a_{1,2} a_{2,3} a_{1,3}$ or 
a cyclic permutation of it.

Now
\begin{eqnarray*}
P(\br a_{1,2} a_{2,3} a_{1,3})&=& P(\br a_{1,2}^{-1} a_{2,3} a_{1,3})
+ z \, P(\br a_{2,3} a_{1,3})\\
&=& P(\br a_{1,2}^{-1} a_{2,3}^{-1} a_{1,3}) + z \, P(\br a_{1,2}^{-1} a_{1,3})
+ z \, P(\br a_{2,3} a_{1,3})\\
&=& P(\br a_{2,3}^{-1}) + z \, P(\br a_{1,2} a_{1,3}) - z^2 \, P(\br a_{1,3})
+ z \, P(\br a_{2,3} a_{1,3}).
\end{eqnarray*} 
By the induction hypothesis the claim follows.
By cyclic permutation the claim also follows for the other two words of length
three. 

Now assume, that $\negat(\br)=1$. 
So, the closed braid can be assumed as $a^{-1} \pi$ where $\pi$ is a
positive braid, with respect to the band generators, and $a$ is one of 
these generators.  
Again the braid is assumed to be of minimal length.

Let $S_+$, $S_-$ and $S_0$ be minimal spanning surfaces for
the closed braids $\br a$, $\br a^{-1}$ and $\br$, where $a$ is a generator
of $B_3$.
Essential to our arguments is a theorem of Scharlemann and 
Thompson \cite{ST:Linkgenus}. 
Two of the three values
\begin{eqnarray}
-\chi(S_+), -\chi(S_-)& \mbox{ and  } & -\chi(S_0)+1 \label{ScharleThompson}
\end{eqnarray}
are the same and the third one is not bigger than the other two.  We 
use the theorem in the following way: In each application of the skein 
relation (\ref{skein}) one of the three terms in 
(\ref{ScharleThompson}) is strictly less than one of the others.  
Hence, the third one has to be equal to the larger one.

Our claim is that the highest degree in $z$ is $l-3$, where $l$ is the minimal 
word length; this equals $-\chi$ of a minimal spanning surface for the braid.
Its coefficient has a non-trivial term in $v^2$. So, Morton's conjecture
follows from this.

If the braid $a^{-1} \pi$ contains a square in $\pi$, then we can reduce the
braid with the help of (\ref{induction 1}). 
Otherwise, we can assume that every positive subword of length three is of 
the form $a_{1,2} a_{2,3} a_{1,3}$ or a cyclic permutation of it. 

If the word length is one then $P(v,z)=v^2 ((1-v^2)/z)^2$.
If the word has length two then an easy computation shows that
$P(v,z)= (v^2-v^4)/z$. 

Now assume, that the word starts with $a_{1,2}^{-1} a_{2,3} a_{1,3}$
We have for some positive word $\pi$
\begin{eqnarray*}
P(a_{1,2}^{-1} a_{2,3} a_{1,3} \pi)&=&
P(a_{1,2}^{-1} a_{2,3}^{-1} a_{1,3} \pi) + z P(a_{1,2}^{-1} a_{1,3} \pi)\\
&=& P(a_{2,3}^{-1} \pi) + z P(a_{1,2}^{-1} a_{1,3} \pi).
\end{eqnarray*}
By the result of Scharlemann and Thompson the word in the 
second term in the last equation has to have minimal braid length and thus
the claim follows by induction.
With the same argument the claim follows if the word starts with
$a_{2,3}^{-1} a_{1,3} a_{12}$ or $a_{1,3}^{-1} a_{1,2} a_{2,3}$.

If the word starts with $a_{1,2}^{-1} a_{1,3} a_{1,2}$ it has to be of 
length at
least $4$, otherwise it wouldn't be minimal. Since we already assumed that
it is square-free, it has to start with
$a_{1,2}^{-1} a_{1,3} a_{1,2} a_{2,3}$,
thus
\begin{eqnarray*}
\lefteqn{P(a_{1,2}^{-1} a_{1,3} a_{1,2} a_{2,3} \pi)=}\\ 
&=&P(a_{1,2}^{-1} a_{1,3}^{-1} a_{1,2} a_{2,3} \pi) + z \, P(a_{2,3} \pi)\\
&=& P(a_{1,2}^{-1} a_{1,3}^{-1} a_{1,2}^{-1} a_{2,3} \pi)
+ z \, P(a_{1,2}^{-1} a_{1,3}^{-1} a_{2,3} \pi) + z \, P(a_{2,3} \pi)\\
&=& P(a_{1,2}^{-2} \pi) +
z \, P(a_{1,2}^{-1} a_{1,3}^{-1} a_{2,3} \pi) + z P(a_{2,3} \pi)\\
&=& P(\pi) + z \, P(a_{1,2}^{-1} \pi)
+ z \, P(a_{1,2}^{-1} a_{1,3} a_{2,3} \pi) - 
z^2\, P(a_{1,2}^{-1} a_{2,3} \pi)  + z \, P(a_{2,3} \pi)\\
&=& P(\pi) + z \, P(a_{1,2}^{-1} \pi)
+ z \, P(a_{1,3} \pi) - 
z^2\, P(a_{1,2}^{-1} a_{2,3} \pi)  + z \, P(a_{2,3} \pi).
\end{eqnarray*}
By induction and using the theorem of Scharlemann and Thompson, it follows 
now that the highest term is coming from
$z^2 \, P(a_{1,2}^{-1} a_{2,3} \pi)$.

If $a_{2,3}^{-1} a_{1,2} a_{2,3}$ or $a_{1,3}^{-1} a_{2,3} a_{1,3}$ are
the initial subwords then the claim follows in the same way.
This completes the proof.
\end{proof}

\subsection{The difference between the Bennequin number and the
Thurston norm}

For a link $L$ let $\max \beta_t(L)$ be the maximal Bennequin number 
among all representatives of $L$ and its mirror image 
and let $S(L)$ be a minimal spanning
surface for $L$. (For definitions see Section  \ref{sec:Benn_inequality})
We give a short argument for the recent result of Kanda \cite{Kanda}:

\begin{proposition}{\rm\cite{Kanda}}\qua
The difference between $\max \beta_t(L)$ and\hfill\break $-\chi(S(L))$ 
can be arbitrarily large.
\end{proposition}

\begin{proof}
For $k$ odd let $\br$ be $\s_2^{-k} \s_1^{k}$. 
Since $\br$ is homogeneous the closure of $\br$ has as maximal
Euler characteristic over all spanning surfaces $3-2 k$.
Obviously $\hat \br$ is equivalent to its mirror image.

It follows immediately from the last section that the difference between the
maximal Bennequin number among any representative of $\hat \br$ and 
the Bennequin number of the special $\br$ is at most $4$.
Thus, $\max \beta_t(L)$ is less than or equal to 
$-3 + 4 =1$. The claim follows.

It can actually be easily shown that if the braid has the 
form $\s_2^{-l} \s_1^{k}$ with positive $k$ and $l \geq 2$ 
then the maximal possible Bennequin number among all representatives
of the closure of the braid is $-l+k-3$.
For this we show that the polynomial $P(0,z)(\s_2^{-l} \s_1^k)$ does not
vanish for $l \geq 2$. Thus, with Corollary \ref{maximal Bennequin number}
we get the desired result.

We already know by the proof of Theorem 
\ref{Morton's conjecture for 3-braids}
that the claim is true
for the closed $2$-braid $\s_1^k$. By Equation 
(\ref{equation circle in HOMFLY}) 
it follows that $$P(0,z)(\id_3 \s_1^k) = 1/z P(0,z)(\s_1^k).$$
Hence, using the defining relations for $P$ we get
\begin{eqnarray*}
P(0,z)(\s_2^{-1} \s_1^k) &=&0\\
P(0,z)(\s_2^{-2} \s_1^k) &=& P(0,z)(\id_3 \s_1^k)-z P(0,z)(\s_2^{-1} \s_1^k)\\
&=& 1/z P(0,z) (\s_1^k)\\
P(0,z)(\s_2^{-l} \s_1^k) &=& P(0,z)(\s_2^{-l+2} \s_1^k) -
z P(0,z)(\s_2^{-l+1} \s_1^k).
\end{eqnarray*}
By induction we can conclude that every $P(0,z)(\s_2^{-l} \s_1^{k})$
is a polynomial $$f(z) P(0,z) (\s_1^k).$$
For $l \geq 3$ odd, all coefficients in $f(z)$ are negative, and
for $l$ even, all coefficients are positive. Thus, cancellation cannot
occur and $f(z)$ is nontrivial for all $l \geq 2$.
\end{proof}

Note, that in contrast to this result, the difference between the
maximal Bennequin number and the minimal degree of the HOMFLY polynomial
in the framing invariant $v$ can be at most $2 n -2$ where $n$ is
the number of strands of the representative.

\Addresses\recd
\end{document}

%% file: agtout.tex
\def\ifplaintex{\expandafter\ifx\csname documentclass\endcsname\relax}

\def\gtp{{\mathsurround=0pt\it $\cal G\mskip-2mu$eometry \&\ 
$\cal T\!\!$opology $\cal P\!$ublications}}  

\def\recd{{\small Received:\qua\receiveddate\ifx\reviseddate\relax
\else\qquad Revised:\qua\reviseddate\fi\par}} 

\def\lognumber#1{\def\thelognumber{#1}}
\def\volumenumber#1{\def\thevolumenumber{#1}}
\def\volumeyear#1{\def\thevolumeyear{#1}}
\def\papernumber#1{\def\thepapernumber{#1}}
\def\pagenumbers#1#2{\def\startpage{#1}\def\finishpage{#2}}
\def\published#1{\def\publishdate{#1}}

\def\received#1{\def\receiveddate{#1}}
\def\revised#1{\def\reviseddate{#1}}
\def\accepted#1{\def\accepteddate{#1}}

\def\asciiaddress#1{\def\theasciiaddress{#1}}
\def\asciiemail#1{\def\theasciiemail{#1}}

\def\asciikeywords#1{\def\theasciikeywords{#1}}

\let\\\par\let\thelognumber\relax\let\thevolumenumber\relax
\let\thepapernumber\relax\let\thevolumeyear\relax\let\startpage\relax
\let\finishpage\relax\let\publishdate\relax\let\receiveddate\relax
\let\reviseddate\relax\let\accepteddate\relax\let\theasciititle\relax
\let\theasciiauthors\relax\let\theasciiaddress\relax
\let\theasciiabstract\relax\let\theasciikeywords\relax

\let\theasciiemail\relax

\ifplaintex
\font\logobig=cmssbx10 scaled 3836
\font\logomed=cmssbx10 scaled 2557
\else
\font\logobig=cmssbx10 scaled 4200
\font\logomed=cmssbx10 scaled 2800
\fi

\long\def\makeagttitle{   
\count0=\startpage
\agt\hfill      
\hbox to 45truept{\vbox to 0pt{\vglue -13truept{\logomed A\kern -.37em{\logobig 
T}\kern -.38em G}\vss}\hss}
\break
{\small Volume \thevolumenumber\ (\thevolumeyear)
\startpage--\finishpage\nl
Published: \publishdate}

\vglue .25truein

{\parskip=0pt\leftskip 0pt plus
1fil\def\\{\par\smallskip}{\Large\bf\thetitle}\par\medskip} \vglue
0.05truein

{\parskip=0pt\leftskip 0pt plus 1fil\def\\{\par}{\sc\theauthors}
\par\medskip}%
 
\vglue 0.03truein 

{\small\leftskip 25truept\rightskip 25truept{\bf Abstract}\stdspace\theabstract

{\bf AMS Classification}\stdspace\theprimaryclass
\ifx\thesecondaryclass\relax\else; \thesecondaryclass\fi\par
{\bf Keywords}\stdspace \thekeywords\par}\vglue 7truept

}   

\ifplaintex
\font\phead=cmsl9 scaled 950
\font\pnum=cmbx10 scaled 913
\font\pfoot=cmsl9 scaled 950
\headline{\vbox to 0pt{\vskip -4.5mm\line{\small\phead\ifnum
\count0=\startpage ISSN 1472-2739 (on-line) 1472-2747 (printed)
\hfill {\pnum\folio}\else\ifodd\count0\def\\{ }%
\ifx\theshorttitle\relax\thetitle\else\theshorttitle\fi\hfill{\pnum\folio}
\else\def\\{ and }{\pnum\folio}\hfill\ifx\theshortauthors\relax\theauthors
\else\theshortauthors\fi\fi\fi}\vss}}
\footline{\vbox to 0pt{\vglue 0mm\line{\small\pfoot\ifnum\count0=\startpage
\copyright\ \gtp\hfill\else
\agt, Volume \thevolumenumber\ (\thevolumeyear)\hfill\fi}\vss}}
\else
\font=cmsl9 scaled 1050
\font\lnum=cmbx10 
\font=cmsl9 scaled 1050
\makeatletter
\def\@oddhead{{\small\lhead\ifnum\count0=\startpage ISSN 1472-2739 
(on-line) 1472-2747 (printed)\hfill {\lnum\number\count0}\else\ifodd\count0
\def\\{ }\ifx\theshorttitle\relax \thetitle \else\theshorttitle\fi\hfill
{\lnum\number\count0}\else\def\\{ and }{\lnum\number\count0}
\hfill\ifx\theshortauthors\relax 
\theauthors\else\theshortauthors\fi\fi\fi}}\def\@evenhead{\@oddhead}
\def\@oddfoot{\small\lfoot\ifnum\count0=\startpage\copyright\ \gtp\hfill\else
\agt, Volume \thevolumenumber\ (\thevolumeyear)\hfill\fi}
\def\@evenfoot{\@oddfoot}
\makeatother
\fi
\let\maketitlepage\makeagttitle
\let\makeshorttitle\maketitlepage
\let\maketitle\maketitlepage

\newwrite\gtoutfile
\long\gdef\makeheadfile{  
{\def\\{, }\def\s{ }
\immediate\openout\gtoutfile head.xxx
\immediate\write\gtoutfile{To: math@arxiv.org}
\immediate\write\gtoutfile{Subject: put OR rep NNNNN:ppppp}
\immediate\write\gtoutfile{--text follows this line--}
\immediate\write\gtoutfile{Proxy-for: \ifx\theasciiauthors\relax
\theauthors\else\theasciiauthors\fi\s<\ifx\theasciiemail\relax\theemail\else\theasciiemail\fi>}
\immediate\write\gtoutfile{\noexpand\\}
\immediate\write\gtoutfile{Authors: \ifx\theasciiauthors\relax
\theauthors\else\theasciiauthors\fi}
{\def\\{ }\immediate\write\gtoutfile{Title: \ifx\theasciititle\relax
\thetitle\else\theasciititle\fi}}
\immediate\write\gtoutfile{Subj-class: GT or SG, GR etc}
\immediate\write\gtoutfile{MSC-class: \theprimaryclass\ifx\thesecondaryclass\relax\else, \thesecondaryclass\fi}
\immediate\write\gtoutfile{Journal-ref: Algebr. Geom. Topol. \thevolumenumber\s
(\thevolumeyear) \startpage-\finishpage}
\immediate\write\gtoutfile{Comments: Published by Algebraic and
Geometric Topology at}
\immediate\write\gtoutfile{\s\s\s  http://www.maths.warwick.ac.uk/agt/AGTVol\thevolumenumber/agt-\thevolumenumber-\thepapernumber.abs.html}
\immediate\write\gtoutfile{\noexpand\\}
\immediate\write\gtoutfile{}
\ifx\theasciiabstract\relax
\immediate\write\gtoutfile{\theabstract}\else
\immediate\write\gtoutfile{\theasciiabstract}\fi
\immediate\write\gtoutfile{}
\immediate\write\gtoutfile{\noexpand\\}
\immediate\write\gtoutfile{}
\immediate\closeout\gtoutfile}}  

\def\maketitlepage{\makeagttitle\makeheadfile}
\let\makeshorttitle\maketitlepage
\let\maketitle\maketitlepage

\def\ifplaintex{\expandafter\ifx\csname documentclass\endcsname\relax}

\def\gtp{{\mathsurround=0pt\it $\cal G\mskip-2mu$eometry \&\ 
$\cal T\!\!$opology $\cal P\!$ublications}}  

\def\recd{{\small Received:\qua\receiveddate\ifx\reviseddate\relax
\else\qquad Revised:\qua\reviseddate\fi\par}} 

\def\lognumber#1{\def\thelognumber{#1}}
\def\volumenumber#1{\def\thevolumenumber{#1}}
\def\volumeyear#1{\def\thevolumeyear{#1}}
\def\papernumber#1{\def\thepapernumber{#1}}
\def\pagenumbers#1#2{\def\startpage{#1}\def\finishpage{#2}}
\def\published#1{\def\publishdate{#1}}

\def\received#1{\def\receiveddate{#1}}
\def\revised#1{\def\reviseddate{#1}}
\def\accepted#1{\def\accepteddate{#1}}

\def\asciiaddress#1{\def\theasciiaddress{#1}}
\def\asciiemail#1{\def\theasciiemail{#1}}

\def\asciikeywords#1{\def\theasciikeywords{#1}}

\let\\\par\let\thelognumber\relax\let\thevolumenumber\relax
\let\thepapernumber\relax\let\thevolumeyear\relax\let\startpage\relax
\let\finishpage\relax\let\publishdate\relax\let\receiveddate\relax
\let\reviseddate\relax\let\accepteddate\relax\let\theasciititle\relax
\let\theasciiauthors\relax\let\theasciiaddress\relax
\let\theasciiabstract\relax\let\theasciikeywords\relax

\let\theasciiemail\relax

\ifplaintex
\font\logobig=cmssbx10 scaled 3836
\font\logomed=cmssbx10 scaled 2557
\else
\font\logobig=cmssbx10 scaled 4200
\font\logomed=cmssbx10 scaled 2800
\fi

\long\def\makeagttitle{   
\count0=\startpage
\agt\hfill      
\hbox to 45truept{\vbox to 0pt{\vglue -13truept{\logomed A\kern -.37em{\logobig 
T}\kern -.38em G}\vss}\hss}
\break
{\small Volume \thevolumenumber\ (\thevolumeyear)
\startpage--\finishpage\nl
Published: \publishdate}

\vglue .25truein

{\parskip=0pt\leftskip 0pt plus
1fil\def\\{\par\smallskip}{\Large\bf\thetitle}\par\medskip} \vglue
0.05truein

{\parskip=0pt\leftskip 0pt plus 1fil\def\\{\par}{\sc\theauthors}
\par\medskip}%
 
\vglue 0.03truein 

{\small\leftskip 25truept\rightskip 25truept{\bf Abstract}\stdspace\theabstract

{\bf AMS Classification}\stdspace\theprimaryclass
\ifx\thesecondaryclass\relax\else; \thesecondaryclass\fi\par
{\bf Keywords}\stdspace \thekeywords\par}\vglue 7truept

}   

\ifplaintex
\font\phead=cmsl9 scaled 950
\font\pnum=cmbx10 scaled 913
\font\pfoot=cmsl9 scaled 950
\headline{\vbox to 0pt{\vskip -4.5mm\line{\small\phead\ifnum
\count0=\startpage ISSN 1472-2739 (on-line) 1472-2747 (printed)
\hfill {\pnum\folio}\else\ifodd\count0\def\\{ }%
\ifx\theshorttitle\relax\thetitle\else\theshorttitle\fi\hfill{\pnum\folio}
\else\def\\{ and }{\pnum\folio}\hfill\ifx\theshortauthors\relax\theauthors
\else\theshortauthors\fi\fi\fi}\vss}}
\footline{\vbox to 0pt{\vglue 0mm\line{\small\pfoot\ifnum\count0=\startpage
\copyright\ \gtp\hfill\else
\agt, Volume \thevolumenumber\ (\thevolumeyear)\hfill\fi}\vss}}
\else
\font=cmsl9 scaled 1050
\font\lnum=cmbx10 
\font=cmsl9 scaled 1050
\makeatletter
\def\@oddhead{{\small\lhead\ifnum\count0=\startpage ISSN 1472-2739 
(on-line) 1472-2747 (printed)\hfill {\lnum\number\count0}\else\ifodd\count0
\def\\{ }\ifx\theshorttitle\relax \thetitle \else\theshorttitle\fi\hfill
{\lnum\number\count0}\else\def\\{ and }{\lnum\number\count0}
\hfill\ifx\theshortauthors\relax 
\theauthors\else\theshortauthors\fi\fi\fi}}\def\@evenhead{\@oddhead}
\def\@oddfoot{\small\lfoot\ifnum\count0=\startpage\copyright\ \gtp\hfill\else
\agt, Volume \thevolumenumber\ (\thevolumeyear)\hfill\fi}
\def\@evenfoot{\@oddfoot}
\makeatother
\fi
\let\maketitlepage\makeagttitle
\let\makeshorttitle\maketitlepage
\let\maketitle\maketitlepage

\newwrite\gtoutfile
\long\gdef\makeheadfile{  
{\def\\{, }\def\s{ }
\immediate\openout\gtoutfile head.xxx
\immediate\write\gtoutfile{To: math@arxiv.org}
\immediate\write\gtoutfile{Subject: put OR rep NNNNN:ppppp}
\immediate\write\gtoutfile{--text follows this line--}
\immediate\write\gtoutfile{Proxy-for: \ifx\theasciiauthors\relax
\theauthors\else\theasciiauthors\fi\s<\ifx\theasciiemail\relax\theemail\else\theasciiemail\fi>}
\immediate\write\gtoutfile{\noexpand\\}
\immediate\write\gtoutfile{Authors: \ifx\theasciiauthors\relax
\theauthors\else\theasciiauthors\fi}
{\def\\{ }\immediate\write\gtoutfile{Title: \ifx\theasciititle\relax
\thetitle\else\theasciititle\fi}}
\immediate\write\gtoutfile{Subj-class: GT or SG, GR etc}
\immediate\write\gtoutfile{MSC-class: \theprimaryclass\ifx\thesecondaryclass\relax\else, \thesecondaryclass\fi}
\immediate\write\gtoutfile{Journal-ref: Algebr. Geom. Topol. \thevolumenumber\s
(\thevolumeyear) \startpage-\finishpage}
\immediate\write\gtoutfile{Comments: Published by Algebraic and
Geometric Topology at}
\immediate\write\gtoutfile{\s\s\s  http://www.maths.warwick.ac.uk/agt/AGTVol\thevolumenumber/agt-\thevolumenumber-\thepapernumber.abs.html}
\immediate\write\gtoutfile{\noexpand\\}
\immediate\write\gtoutfile{}
\ifx\theasciiabstract\relax
\immediate\write\gtoutfile{\theabstract}\else
\immediate\write\gtoutfile{\theasciiabstract}\fi
\immediate\write\gtoutfile{}
\immediate\write\gtoutfile{\noexpand\\}
\immediate\write\gtoutfile{}
\immediate\closeout\gtoutfile}}  

\def\maketitlepage{\makeagttitle\makeheadfile}
\let\makeshorttitle\maketitlepage
\let\maketitle\maketitlepage